\documentclass[a4paper,12pt]{amsart}

\usepackage{amssymb,amscd,amsthm,mathrsfs}
\usepackage[ansinew]{inputenc}
\usepackage[centertags]{amsmath}

\newcommand{\en}{\subset}

\newcommand{\Z}{\mbox{$\mathbb{Z}$}}

\newcommand{\Res}{\mbox{$\mathcal{R}$}}
\newcommand{\T}{\mbox{$\mathbb{T}$}}
\newtheorem{teo}{Theorem}
\newtheorem{conj}{Conjecture}
\newtheorem{cor}{Corollary}
\newtheorem{lema}{Lemma}

\newcommand{\bi}{\begin{itemize}}
\newcommand{\ei}{\end{itemize}}

\theoremstyle{definition}

\theoremstyle{remark}

\newcommand{\eps}{\varepsilon}
\newcommand{\la}{\lambda}

\newcommand{\dem}[1]{\vspace{.05in}{\sc\noindent Proof #1.}}
\newcommand{\lqqd}{\par\hfill {$\Box$} \vspace*{.1in}}

\newcommand{\diff}[1]{\mathcal{D}iff^{#1}(M)}
\newcommand{\U}{\mathcal{U}}

\author[R. Potrie]{Rafael Potrie}
\address{CMAT, Facultad de Ciencias, Universidad de la Rep\'ublica, Uruguay}
\email{rpotrie@cmat.edu.uy}
\author[M. Sambarino]{Martin Sambarino}
\address{CMAT, Facultad de Ciencias, Universidad de la Rep\'ublica, Uruguay}
\email{samba@cmat.edu.uy}

\title{Codimension one generic homoclinic classes with interior}

\begin{document}

\maketitle

\begin{abstract}
We study $C^1$-generic diffeomorphisms with a homoclinic class with non
empty interior and in particular those admitting a codimension one
dominated splitting. We prove that if in the finest dominated
splitting the extreme subbundles are one dimensional then the
diffeomorphism is partially hyperbolic and from this we deduce that
the diffeomorphism is transitive.
\end{abstract}

\section{Introduction}

It is a main problem in generic dynamics to understand the structure of homoclinic classes, this has became most important after some results in \cite{BC} which raised the interest in the study of chain recurrence classes and in particular homoclinic classes which are, generically, the chain recurrence classes containing periodic points. This results, as most of the results (including the ones here presented) are in the $C^1$ category, very little is known about $C^r$ generic diffeomorphisms with $r>1$ (see \cite{Pujals}).

A lot is known in the case of a isolated homoclinic class $H(p,f)$
of a hyperbolic periodic point $p$ (i.e., it is maximal invariant in
a neighborhood) of a \textit{generic diffeomorphism $f$}, see for
instance  \cite{beyondhip}, chapter 10 (we remark here  that the
genericity of $f$ implies that the continuation $H(p_g,g)$ is also
isolated for $g$ in a neighborhood of $f$). The key point of knowing
that the class is isolated is that, after perturbation, orbits that
remains in a neighborhood must belong to the continuation of the
homoclinic class $H(p_g,g).$ However, if one does not know
\textit{apriori} that the homoclinic class is isolated (in other
words, the class might be wild, i.e., non isolated) only very sparse
results have been obtained (see for example \cite{ABCDW}
 where they prove that generic homoclinic classes are index complete). The main difficulty is to
 overcome the fact that after perturbations one cannot ensure that the perturbed points remain in the class.

For example, it is not known whether the non-wandering set of a $C^1$-generic diffeomorphism may have nonempty interior and not coincide with the whole manifold. We treat this problem in this paper solving it in some particular cases and also giving some results which may help to obtain a general solution. We remark that it is not difficult to construct examples of diffeomorphisms such that its non-wandering set has non empty interior and doesn't coincide with the whole manifold.

\subsection{Definitions and statement of results}

\smallskip

Let $M$ be a compact connected boundaryless manifold of dimension
$d$ and let $\diff 1$ be the set of diffeomorphisms of $M$ endowed
with the $C^1$ topology. We shall say that a property (or a diffeomorphism) is \emph{generic} if
and only if there exists a residual ($G_{\delta}$-dense) set $\Res$ of $\diff 1$ for
which for every $f \in \Res$ satisfies that property.

The main result of this paper concerns the following conjecture of
\cite{ABD}, we shall denote $\Omega(f)$ to the nonwandering set of $f$ :

\begin{conj}
 There is a residual set $\Res \en \diff 1$ such that if $f\in \Res$ and $int(\Omega(f))\neq \emptyset$ then $f$ is transitive.
\end{conj}

For a hyperbolic periodic point $p \in M$ of some diffeomorphism $f$
we denote its \emph{homoclinic class} by $H(p,f)$, defined as the
closure of the transversal intersections between the stable and
unstable manifolds of the orbit of $p$. For generic diffeomorphisms, if the nonwandering set has nonempty interior, then there is a homoclinic class with nonempty interior (see \cite{ABD} or \cite{BC}). So, the conjecture is reduced to the study of homoclinic classes with nonempty interior.

Some progress has been made towards the proof of this conjecture
(see \cite{ABD} and \cite{ABCD}), in particular, it has been proved
in \cite{ABD} that isolated homoclinic classes as well as homoclinic
classes admitting a strong partially hyperbolic splitting (we shall define this concept later) verify the
conjecture. Also, they proved that a homoclinic class with non empty
interior must admit a dominated splitting (see Theorem 8 in
\cite{ABD}). In \cite{ABCD} the conjecture was proved for surface
diffeomorphisms.

In \cite{ABD} the question about whether within the finest dominated
splitting the extremes subbundles should be volume hyperbolic was
posed. We give a positive answer when the class admits codimension
one dominated splitting. This gives also new situations where the
above conjecture holds and weren't known.

Let us recall the definition of dominated splitting: a compact set
$H$ invariant under a diffeomorphism $f$ admits dominated splitting
if the tangent bundle over $H$ splits into two $Df$ invariant
subbundles $T_HM=E\oplus F$ such that there exist $C>0$ and
$0<\lambda<1$ such that for all $x\in H:$
$$\|Df^n_{/E(x)}\|\|Df^{-n}_{/F(f^n(x))}\|\le C\lambda^n$$

\noindent we say in this case that $F$ dominates $E$.

Let us remark that Gourmelon (\cite{GourmelonAdaptada}) proved that
there  always exists an adapted metric for which $C=1.$

We shall say that a bundle $F$ is uniformly expanding (contracting) if there exists $n_0<0$ ($n_0>0$) such that $\|Df^{n_0}_{/F(x)}\|< 1/2$ $\forall x \in H$.

The main theorem of this paper is the following

\begin{teo}\label{maintheorem}
 Let $f$ be a generic diffeomorphism with a homoclinic class $H$ with non empty interior
 and
 admitting a codimension one dominated splitting $T_H M=E \oplus F$
 where $dim(F)=1$. Then, the bundle $F$ is uniformly expanding for $f$.
\end{teo}

As a consequence of our main theorem we get the following easy corollaries.

Recall that a compact invariant set $H$ is strongly partially hyperbolic if it admits a three ways dominated splitting $T_HM= E^s \oplus E^c \oplus E^u$ (that is, $E^s\oplus E^c$ is dominated by $E^u$ and $E^s$ is dominated by $E^c\oplus E^u$), where $E^s$ is non trivial and uniformly contracting and $E^u$ is non trivial and uniformly expanding.

\begin{cor}\label{Cor1} Let $H$ be a homoclinic class with non empty interior
for a generic diffeomorphism $f$ such that $T_HM= E^1 \oplus E^2
\oplus E^3$ is a dominated splitting for $f$ and
$dim(E^1)=dim(E^3)=1$. Then, $H$ is strongly partially hyperbolic and $H=M$.
\end{cor}

\dem{} The class should be strongly partially hyperbolic because of
the previous theorem (applied to $f$ and to $f^{-1}$).  Corollary 1
of \cite{ABD} (page 185)  implies that $H=M$. \lqqd

We say that a homoclinic class $H$ is far from tangencies if there
is a neighborhood of $f$ such that there are no homoclinic
tangencies associated to periodic points in the continuation of $H$. The tangencies are of index $i$ if they are associated to a periodic point of index $i$, that is, its stable manifold has dimension $i$.
We  get the following result following the generalization of the results of \cite{W} in \cite{Gou2} (see also
\cite{ABCDW}):

\begin{cor} Let $H$ be a homoclinic class with non empty interior
for a generic diffeomorphism $f$ such that $H$ is far from tangencies of index $1$ and $n-1$ and has index $1$ and
$n-1$ periodic points. Then, $H=M$.
\end{cor}

\dem{} Since the class is far from tangencies, and the classes for
generic diffeomorphisms either  coincide or are disjoint (see
\cite{BC}) we have that using \cite{Gou2} the class must admit a
dominated splitting with one dimensional extremal subbundles (see
also \cite{ABCDW} Corollary 3), thus, by using Corollary \ref{Cor1}
we get the result.

\lqqd

In fact, the previous corollary can be compared to a corollary of a new result of Yang (\cite{Y}) on Lyapunov stable homoclinic classes far from tangencies. Yang's result (Theorem 3 in \cite{Y}) implies that a generic homoclinic class with nonempty interior and far away from (any) tangencies must be strongly partially hyperbolic or contained in the closure of the set of sinks and sources, and thus (using the results of \cite{ABD}) the whole manifold.

Incidentally, we also give a new proof in the two dimensional case:

\begin{cor}
Let $f$ be a generic surface diffeomorphism having a homoclinic
class with nonempty interior. Then $f$ is conjugated to a linear Anosov diffeomorphism in $\T^2$.
\end{cor}

\dem{} Since the class must admit dominated splitting (Theorem 8 of \cite{ABD}), this should be into 2 one dimensional subbundles. So, the class must be hyperbolic and thus, since the conjecture holds for hyperbolic homoclinic classes $f$ is Anosov (the rest follows from classical theory of Anosov diffeomorphisms).
\lqqd

\subsection{Idea of the proof}

\smallskip

The idea of the proof is the following.

First we prove that if the homoclinic class has interior, the
periodic points in the class (which are all saddles) should have
eigenvalues (in the $F$ direction) exponentially (with the period)
far from $1$. Otherwise we manage to obtain a sink or a source
inside the interior of the class and thus contradicting the fact
that the interior of the homoclinic class for generic
diffeomorphisms is, roughly speaking, robust (Theorem 4 of
\cite{ABD}).

Then, using the previous fact and some results of \cite{pablito} and
\cite{integrability} we manage to prove that  the center manifolds
integrating  a one dimensional extreme subbundle  should have nice
dynamical properties. For this we also use the connecting lemma for
pseudo orbits of \cite{BC}.

Finally, in the event that the extreme subbundle is not hyperbolic,
we manage to obtain (using dynamical properties and a Lemma of Liao)
periodic points near the class with bad contraction or expansion in
those extreme subbundles. Using Lyapunov stability of the homoclinic
class (which is generic, see \cite{ABD} and \cite{CMP}) we  ensure
that the periodic points we found belong to the class and thus reach
a contradiction.

\medskip

\textit{Acknowledgements:} We would like to thank Andr\'es Sambarino for
motivating us to study this problem. Also, we thank Sylvain Crovisier for reading a draft, explaining us Yang's results and \cite{C}.
Finally, we would like thank the referee for many valuable comments that helped improve the presentation of this paper.

\section{Preliminary results}

In this section we shall state some results we are going to use in
the proof of the main theorem. It can be skipped and used as
reference when the results are used.

Some generic properties of diffeomorphisms are contained in the
following Theorem (see \cite{ABD} and references therein):

\begin{teo}\label{propiedadesgenericas}
 There exists a residual subset $\Res$ of $\diff 1$ such that if $f \in \Res$
\begin{itemize}
 \item[a1)] $f$ is Kupka Smale (that is, all its periodic
  points are hyperbolic and their invariant manifolds intersect transversally).
 \item[a2)] The periodic points of $f$ are dense in the chain
 recurrent set of $f$(\footnote{The chain recurrent set is
 the set of points $x$ satisfying that for every $\eps>0$
 there exist an $\eps$-pseudo orbit form $x$ to $x$, that is,
 there exist points $x=x_0, x_1, \ldots x_k=x$, $k>0$ such that
 $d(f(x_i),x_{i+1})< \eps$.}). Moreover, homoclinic classes
 coincide with those chain recurrent classes which contain periodic points.
 \item[a3)] Every homoclinic class with non empty interior of
 $f$ is Lyapunov stable for $f$ and $f^{-1}$ (\footnote{Lyapunov
 stability of $\Lambda$ means that $\forall U$ neighborhood of
 $\Lambda$ there is $V\en U$ neighborhood of $\Lambda$ such that
 $f^n(V) \en U$ $\forall n\geq 0$.}). This implies that the stable and unstable set of any point in the class is contained in the class.
 \item[a4)] For every periodic point $p$ of $f$,
 $H(p,f) =\overline{W^s(p)} \cap \overline{W^u(p)}$.
  \item[a5)] Given a homoclinic class $H$ of a periodic point $p$, if $U$ is an
 open set such that $\overline{U} \en int(H)$
then there exists $\U$ neighborhood of $f$ such that for every $g
\in \U \cap \Res$, $U \en H(p_g,g)$ is satisfied (where $p_g$ is the
continuation of $p$ for $g$).
  \item[a6)] Homoclinic classes vary continuously with the Hausdorff distance with respect to $f$. This means, that given $p\in H$ a periodic point and $\eps>0$, there exists $\U$ a neighborhood of $f$ such that for every $g\in \U$, the homoclinic class of the continuation of $p$ lies within less than $\eps$ from $H$ in the Hausdorff distance.
\end{itemize}
\end{teo}

To obtain dynamical properties of the center manifolds we shall use
the following results from \cite{pablito} and \cite{integrability}.
First recall that if $T_H M= E \oplus F$ is a dominated
splitting then, Theorem 5.5 of \cite{invariant} gives us a local
$f-$invariant manifolds $W^F_{\eps}$ tangent to $F$.

Local $f-$invariance means that $\forall \eps>0$ there exists
$\delta>0$ such that $f^{-1}(W^F_\delta(x)) \en W^F_\eps(f(x))$. Taking
$f^{-1}$ we have an analog for $E$.

\begin{teo}[Main Theorem of \cite{pablito}]\label{pablito}
Let $\Lambda$ a compact invariant set of a generic diffeomorphism
$f$ admitting a codimension one dominated splitting $T_\Lambda M=
E \oplus F$ with $dim(F)=1$. Assume that
$\overline{Per(f_{/\Lambda})}=\Lambda.$ Then, $\forall x \in \Lambda$ and
$\forall \eps>0$ there exists $\delta>0$ such that

\[ f^{-n}(W^F_{\delta}(x)) \en W^F_{\eps}(f^{-n}(x))  \ \ \forall n \geq
0 \]

In particular, $W^F_{\delta}(x) \en \{ y \in M \ : \ d(f^{n}(x),
f^n(y))\leq \eps \}$.
\end{teo}

Remark 2.8 of \cite{C} gives a similar result that is enough in our context.

If there is a dominated splitting for $H$ of the form $T_H M=
E\oplus F$, then, there exists $V$, a neighborhood of $H$ such that
if a point $z$ satisfies that $f^n(z) \in V$ $\forall n \in \Z$ then
we can define the splitting for $z$ and it will be dominated (see
\cite{beyondhip}). Such a neighborhood will be called adapted.

If $I$ is an interval, we denote by $\omega(I) = \bigcup_{x\in I}
\omega(x)$, and by $W^{ss}_\eps(I)= \bigcup_{x\in I} W^{ss}_\eps
(x)$ its strong stable manifold. Also $\ell(I)$ denote its length.
We shall state the following result which is an immediate Corollary
of Theorem 3.1 of \cite{integrability} for generic dynamics.

\begin{teo}[\cite{integrability}]\label{pelado} Let $f\in \diff 1$ a
generic diffeomorphism and $\Lambda$ compact invariant set admitting
a codimension one dominated splitting  $T_\Lambda M= E  \oplus
F$ (where $dim(F)=1$). Then, there exists $\delta_0$ such that
if $I$ is an interval integrating the subbundle $F$ satisfying
$\ell(f^n(I))< \delta <\delta_0$ $\forall n\geq 0$ and that its
orbit remains in an adapted neighborhood $V$ of $\Lambda$, then, only
one of the following holds:

\begin{enumerate}
\item $\omega(I)$ is contained in the set of periodic points
of $f$ restricted to $V$ of $\Lambda$ one of which is an attractor.
\item $I$ is wandering (that is, $W^{ss}_\eps(f^n(I)) \cap W^{ss}_\eps(f^m(I))=
\emptyset$ for all $n\neq m$). This implies that $\ell(f^n(I)) \to
0$ as $|n|\to \infty$.
\end{enumerate}
\end{teo}

Other result we shall use is the following well known Lemma of
Franks:

\begin{teo}[Frank's Lemma \cite{frankslema}]\label{frankslema}
Let $f \in Diff^1(M)$. Given $\mathcal U (f)$, a $C^1$ neighborhood
of $f$, then $\exists \;\mathcal U_0 (f)$ and $\eps>0$ with the
following property: if $g \in \mathcal U_0 (f)$, $\theta= \{
x_1,\ldots, x_m\}$ and

\[ L:\bigoplus_{x_i \in \theta} T_{x_i}M \to \bigoplus_{x_i \in
\theta}T_{g(x_i)} M  \ \ \text{such that} \ \ \left\|L-
Dg|_{\bigoplus T_{x_i}M} \right\|< \eps \]

\noindent are given, then there exists $\tilde{g} \in \mathcal U(f)$
such that $D\tilde{g}_{x_i} = L|_{T_{x_i}M}$. Moreover if $R$ is a
compact set disjoint from $\theta$ we can consider $\tilde{g}=g$ in
$R \cup \theta$.
\end{teo}

Finally we state the following Lemma of Liao. A proof can be found
(with the same notation) in \cite{wen}. We shall state the result in
the particular case of index one dominated splitting with an adapted
metric (which always exist because of \cite{GourmelonAdaptada}), but
it holds in a wider context. Recall also that for linear maps $A_i$
in one dimensional spaces it holds that $\prod_i \|A_i\|= \|\prod_i
A_i \|$.

\begin{lema}[Liao \cite{liao}]\label{lemaliao} Let $\Lambda$ be a compact
invariant set of $f$ with dominated splitting $T_H M= E \oplus
F$ such that $\|Df_{/E(x)}\|\|Df^{-1}_{/F(x)}\|< \gamma$
$\forall x  \in \Lambda$ and $dim (F)=1$. Assume that
\begin{enumerate}
 \item There is a point $b\in \Lambda$ such that $\|Df^{-n}|_{F(b)}\| \geq 1$
 $\forall n\geq 0$.
 \item There exists $\gamma<\gamma_1<\gamma_2<1$ such that given $x\in \Lambda$
  satisfying $\|Df^{-n}_{/F(x)}\|\geq \gamma_2^n$ $\forall n\geq 0$ we have that
   there is $y \in \omega(x)$ satisfying $\|Df^{-n}_{/F(y)}\|\leq \gamma_1^n$
    $\forall n\geq 0$.
\end{enumerate}

Then, for any $\gamma_2<\gamma_3<\gamma_4<1$ and any neighborhood
$U$ of $\Lambda$ there exists a periodic point $p$ of $f$ whose
orbit lies in $U$, is of the same index as the dominated splitting
and satisfies $\|Df^{-n}_{/F(p)}\| < \gamma_4^n$ $\forall n\geq 0$
and $\|Df^{-n}_{/F(p)}\| \geq \gamma_3^n$ $\forall n \geq 0$.

\end{lema}

\section{Proof of the main theorem}

For $p\in Per(f)$, $\pi(p)$ denotes the period of $p$.

\begin{lema}\label{normadiferencial} Let $H$ be a homoclinic class
with interior of a generic diffeomorphism $f$ such that $T_HM=E\oplus F$ is a dominated splitting with $dim F=1$ . Then, there exists $\lambda<1$ such that for all $p \in Per(f|_H)$ the following holds:

$$\|Df_{/F(p)}^{-\pi(p)}\|\le \lambda^{\pi(p)}$$

\end{lema}

\dem{}
Arguing by contradiction assume that the conclusion does not hold, that is,
 for every $\lambda<1$ there exists $p \in
Per(f_{/H})$  such that $\|Df_{/F(p)}^{-\pi(p)}\|\ge
\lambda^{\pi(p)}$ which is equivalent to
$\|Df_{/F(p)}^{\pi(p)}\|\le \la^{-\pi(p)}$ since $F$ is one
dimensional. When the class is isolated this is enough since one can perturb the orbit in order to create a sink contradicting the isolation. Here, since the class may be wild, the creation of the sink represents no contradiction, so we must use the persistence of the interior given by generic property $a5)$ of Theorem \ref{propiedadesgenericas} and create a sink there.

Let $U$ be an open set such that $\overline{U}\en int(H)$. Since $f$
is generic, property $a5)$ of Theorem \ref{propiedadesgenericas}
ensure us the existence of a neighborhood $\U$ of $f$ such that for
every $g$ in a residual subset of $\U$ we have $U \en H_g$ ($H_g$ is
the continuation of $H$ for $g$, from $a5)$ of Theorem \ref{propiedadesgenericas} this continuation makes sense since it will be the only class containing $U$).

Frank's Lemma implies the existence of $\eps>0$ such that if we fix
an arbitrary finite set of points, we can perturb the diffeomorphism
as near as we want of those points obtaining a new diffeomorphism
with arbitrary derivatives ($\eps-$close to the originals) in those points and such that the diffeomorphism lies inside
$\U$.

Let us fix $1>\lambda> 1-\eps/2$  and let $p \in Per(f_{/H})$ as
before.
 Since $f$ is generic, the periodic points of
the same index as $p$ are dense in $H$ so, we can choose $q\in U
\cap Per(f)$ homoclinically related to $p$.

Let $x \in W^s(p)\cap W^u(q)$ and $y \in W^s(q) \cap W^u(p)$, we get
that the set $\Lambda= \mathcal O (p) \cup  \mathcal O (q) \cup
\mathcal O (x) \cup   \mathcal O (y)$ hyperbolic.

Consider  the following periodic pseudo orbit contained in
$\Lambda$, $$\{...,p, f(p),...,f^{N\pi(p)-1}(p), f^{-n_0}(y), \ldots
,f^{n_0}(y), f^{-n_0}(x), \ldots, f^{n_0}(x),p,...\}$$ which we
shall denote as $\wp^N$. Clearly, given $\beta>0$ there exists $n_0$
such that  $\wp^N$ is a $\beta$-pseudo orbit. At the same time, if
we choose $N$ large enough we  obtain a pseudo orbit which stays
near $p$ much longer than of $q$ and then inherit the behavior of
the derivative of $p$ rather than that of $q$.

The shadowing lemma for hyperbolic sets (see \cite{shub}) implies
that for every $\alpha>0$ there exists $\beta$ such that every
closed $\beta$-pseudo orbit is $\alpha$-shadowed by a periodic
point. So, let us choose $\alpha$ in such a way that the following
conditions are satisfied:

\begin{itemize}
\item[(a)] $B_{2\alpha}(q) \en U$.
\item[(b)] If $d(z,w)<\alpha$ and $x,y$ are
in an adapted neighborhood of $H$ then,
$$\frac{\|Df_{/F(z)}\|}{\|Df_{/F(w)}\|}< 1+c$$ ($c$ verifies
$(1+c)(1-\frac{\eps}{2})^{-1}< 1+\eps$).
\end{itemize}

Let $\beta<\alpha$ be given from  the Shadowing Lemma for that
$\alpha$ and let $n_0$ be such that $\wp^N$ is a $\beta$-pseudo
orbit. Therefore there exists a periodic orbit $r$ of period
$\pi(r)=N\pi(p)+4n_0$ which $\alpha$- shadows $\wp^N.$
Therefore, setting $k=\sup_{x\in M} \|Df_x\|$, we have
\begin{eqnarray*}
\|Df_{/F(r)}^{N\pi(p)+4n_0}\|&\le
&k^{4n_0}(1+c)^{N\pi(p)}\|Df_{/F(p)}^{\pi(p)}\|^N\\
&\le &
k^{4n_0}\left((1+c)(1-\frac{\eps}{2})^{-1}\right)^{N\pi(p)}<(1+\eps)^{\pi(r)}
\end{eqnarray*}
where the last inequality holds provided  $N$ is large enough.
Notice that the orbit of $r$ passes through $U.$ On the other hand,
by domination, we have that
$\|Df_{/E(r)}^{\pi(r)}\|<\|Df_{/F(r)}^{\pi(r)}\|.$ Since $E$
and $F$ are invariant we conclude that any eigenvalue of
$Df_r^{\pi(r)}$ is less than $(1+\eps)^{\pi(r)}.$

Now, if we compose in the orbit of $r$ its derivatives with
homoteties of value $(1+\eps)^{-1}$  we obtain, by using Frank's
Lemma, a diffeomorphism $g$ so that all the eigenvalues associated
to the periodic orbit $r$ are less than $1$, that is, $r$ is a
periodic attractor (sink). This contradicts the generic assumption,
since the sink is persistent, so every residual $\Res \in \U$ will
have diffeomorphisms with a sink near $r$, thus contained in $U$,
and thus contradicting that the interior is persistent.

\lqqd

\begin{lema}\label{variedadesdinamicamentedef}
Let $H$ be a homoclinic class with non empty interior for a generic
diffeomorphism $f$ such that $T_HM= E \oplus F$ is a dominated
splitting with $dim(F)=1$. Then, there exists $\eps_0$ such
that for all $\eps\le \eps_0$ there exists $\delta$ such that
$\forall x \in H$, $$W^F_{\delta}(x) \en W^{uu}_{\eps}(x):= \{ y \in M
: \ d(f^{-n}(x), f^{-n}(y))\leq \eps \ ; \ d(f^{-n}(x), f^{-n}(y))\to 0\}$$.
\end{lema}

\dem{} First we shall prove the Lemma for periodic points and then,
using this fact prove the general statement. Let $\eps_0>0$ such
that $B_{\eps_0} (H)$ is contained in the adapted neighborhood of
$H$ and such that if $d(x,y)<\eps_0$ then

\[\frac{\|Df^{-1}_{/F(x)}\|}{\|Df^{-1}_{/F(y)}\|} < \lambda^{-1} \]

\noindent where $\lambda$ is given by Lemma \ref{normadiferencial}.
Let $\eps\le \eps_0$ and let $\delta>0$ from Theorem \ref{pablito}
corresponding to this $\eps.$

Let $p \in Per(f_{/H})$ for which there is $y \in W^F_\delta(p)$ such
that $d(f^{-n}(y), f^{-n}(p)) \nrightarrow 0$. Since $W^F_\delta(p)$
is one dimensional, $W^F_\delta(p) \backslash \{p\}$ is a disjoint
union of two intervals. Denote $I_\delta$ the connected component of
 $W^F_\delta(p) \backslash \{p\}$ that contains $y.$ By Theorem
 \ref{pablito} we have either $f^{2\pi(p)}(I_\delta)\subset
 I_\delta$ or $f^{2\pi(p)}(I_\delta)\supset
 I_\delta.$ In any event, since $y\in I_\delta$ we conclude that
 there exits a point $z_0\in W_\eps^F(p)$ fixed under $f^{2\pi(p)}$
 and such that $\|Df_{/F(z_0)}^{2\pi(p)}\|\le 1.$

This contradicts the previous Lemma, since by the way  $\eps$ was
chosen we get (since  we know that $d(f^i(p), f^i(z_0))<\eps$ for
all $i$) that

\[ \|Df^{2\pi(p)}_{/F(p)}\| = \prod_{i=0}^{2\pi{p}-1} \|Df_{/F(f^i(p))}\|
< \lambda^{-2\pi(p)} \prod_{i=0}^{2\pi{p}-1} \|Df_{/F(f^i(z_0))}\|
=
\]\[= \lambda^{-2\pi(p)} \| Df^{2\pi(p)}|_{F(z_0)}\| < \lambda^{-2\pi(p)} \]

Now, lets prove the general statement. Let us suppose that for every
$\eps>0$ there exist $x\in H$ and a small arc $I\en W^F_\delta(x)$
containing $x$ such that $\ell(f^{-n}(I))\nrightarrow 0$. We know,
because of Theorem \ref{pablito} that $\ell(f^{-n}(I))\leq \eps$,
then, taking $n_j\to +\infty$ such that $\gamma\leq
\ell(f^{-n_j}(I))\leq \eps$ and taking limits, we obtain an arc $J$
integrating $F$ such that $\ell(f^{n}(J))\leq \eps$ $\forall n
\in \Z$ and containing a point $z \in J\cap H$ (a limit point of
$f^{-n_j}(x)$).

Now, we shall use Theorem \ref{pelado} to reach a contradiction. It
is not difficult to discard the first possibility in the Theorem
because it will contradict what we have proved for periodic points.

On the other hand, if $J$ is wandering, we know that it can not be
accumulated by periodic points. Since $f$ is generic, we reach a
contradiction if we prove that the points in $J$ are chain recurrent
(see property $a2)$ of Theorem \ref{propiedadesgenericas}). Theorem
\ref{pelado}, implies that, $\ell(f^n(J)) \to 0$ ($|n| \to
+\infty$), then, since $z \in H\cap J$, if we fix $\epsilon$, and $y
\in J$, then, for some future iterate  $k_1$ and a past one $-k_2$,
we know that $f^{k_1}(y)$ is $\epsilon$-near of $f^{k_1}(z)$ and
$f^{-k_2}(y)$ is $\epsilon$-near  $f^{-k_2}(z)$. Since homoclinic
classes are chain recurrent classes, there is an $\epsilon$ pseudo
orbit from $f^{k_1}(z)$ to $f^{-k_2}(z)$ and then, $y$ is chain
recurrent, a contradiction.

\lqqd

\begin{cor}  Let $H$ be a homoclinic class with non empty interior
for a generic diffeomorphism $f$ such that $T_HM= E \oplus F$ is
a dominated splitting with $dim(F)=1$. Then, $F$ is
uniquely integrable. \end{cor}

\dem{} It follows from the fact that the center stable manifold is
dynamically defined (see  \cite{invariant}). \lqqd

Uniqueness of the center manifolds imply that one can know that if you have a point $y\in W^F_\delta (x) \cap H$ then there exists $\gamma<\delta$ such that $W^F_\gamma (y)\en W^F_\delta (x)$.

\begin{cor}\label{vaacero} Let $H=H(p,f)$ be a homoclinic class with non empty interior
for a generic diffeomorphism $f$ such that $T_HM=E \oplus F$ is
a dominated splitting with $dim(F)=1$. Then, for all $L>0$
and $l>0$ there exists $n_0$ such that if $I$ is a compact arc
integrating $F$ whose length is smaller than $L$, then
$\ell(f^{-n}(I))< l $ $\forall n>n_0$.
\end{cor}

\dem{} It is easy to see that every compact arc integrating $F$
should have its iterates of length going to zero in the past because
of Theorem \ref{variedadesdinamicamentedef} (it is enough to
consider a finite covering of $I$ where the Theorem applies).

Lets suppose then that there exists $L$ and $l$ such that for every
$j>0$ there is an arc $I_j$ integrating $F$ of length smaller than
$L$ and $n_j>j$ such that $\ell(f^{-n_j}(I_j) \geq l$. We can
suppose without loss of generality that $\ell(I_j) \in (L/2,L)$.

Also, we can assume (maybe considering subsequences) that $I_j$
converges uniformly to an arc $J$ integrating $F$ and verifying
$L/2\leq \ell(J) \leq L$.

Since the length of $J$ is finite and it integrates $F$ we know
that $\ell(f^{-n}(J)) \to 0$ with $n\to +\infty$.

Let $\eps= l/2$ and $\delta$ given by Theorem \cite{pablito} which
ensures that $W^F_\delta(x) \en W^u_\eps(x)$ $\forall x$.

Let $n_0$ such that $\forall n \geq n_0$ we have $\ell(f^{-n}(J)) <
\delta/4$. Let also be $\gamma$ small enough such that if $x\in
B_{\gamma}(J)$ then $d(f^{-k}(x),f^{-k}(J))< \delta/4$ $\forall
0\leq k \leq n_0$.

Now, if we consider $j$ large enough (in particular $j>n_0$) such
that $I_j \en B_\gamma(J)$ we obtain $\ell(f^{-n_0}(I_j)) < \delta$
and so $\ell(f^{-n}(I_j)) < \eps < l$ $\forall n \geq n_0$, so,
$n_j< n_0$ which is a contradiction.

\lqqd

We are ready to give the proof of our main theorem:

\begin{teo} Let $H$ be a homoclinic class with non empty interior
for a generic diffeomorphism $f$ such that $T_HM= E \oplus F$ is
a dominated splitting with $dim(F)=1$. Then, $F$ is
uniformly expanding 
\end{teo}

\dem{} Because of the existence of an adapted norm for the
dominated splitting (see \cite{GourmelonAdaptada}) we can assume
that $\|Df|_{E(x)}\|\|Df^{-1}|_{F(f(x))}\| < \gamma$ (for the
sake of simplicity).

Suppose the theorem is not true. Thus,  for every $0<\nu<1$ there
exists some $x\in H$ such that $\|Df^{-n}_{/F(x)}\| \geq \nu$,
$\forall n\geq 0$ (otherwise for every $x$ there would be some
$n_0(x)$ which would be the first one for which $\|Df^{-n}_{/F(x)}\|<
\nu$ and by compactness, the $n_0(x)$ are uniformly bounded, then $F$
would be hyperbolic). If we choose points $x_m$ satisfying
$\|Df^{-n}_{/F(x)}\|\geq 1 - 1/m$ $\forall n \geq 0$, a limit
point $x$ will satisfy $\|Df^{-n}_{/F(x)}\|\geq 1$ $\forall n \geq 0$.

First of all, we consider the case where we cannot use the Shifting
Lemma of Liao (Lemma \ref{lemaliao}). It is not difficult to see
that this implies (using Pliss' Lemma, see also \cite{wen}) that
$\forall \gamma <\gamma_1< \gamma_2 < 1$, there exists $x\in H$ such that 

\[ \|Df^{-n}_{/F(x)}\| \geq \gamma_2^n \qquad \forall n \geq 0 \]

\noindent but, $\forall y \in \omega(x)$ we have that

\[ \|Df^{-n}_{/F(y)}\| \geq \gamma_1^n \qquad \forall n \geq 0 \]

So, if we work in $\omega(x)$ which is a closed invariant set, we
have that the subbundle $E$ will be hyperbolic since the dominated
splitting implies that $\forall z \in \omega(x)$

\[ \|Df_{/E(z)}\| < \frac{\gamma}{\|Df^{-1}_{/F(f(z))}\|}<
\frac{\gamma}{\gamma_1} < 1 \]

This implies that, since we have dynamical properties for the
manifolds integrating the subbundle $F$, that we can  shadow
recurrent orbits. Indeed, if we have a recurrent point $y \in
\omega(x)$, for every small $\eps$ (in particular, such that the
stable and unstable manifolds of $y$ are well defined) we can
consider $n$ large enough so that $d(f^n(y),y) \leq \eps/3$,
$f^n(W^E_\eps (y)) \en W^E_{\eps/3}(f^n(y))$ and
$f^{-n}(W^F_{\eps}(f^n(y))) \en W^F_{\eps/3}(y)$ which gives us
(using classical arguments) a periodic point $p$ of $f$ which
verifies that has period $n$ and remains $\eps$-close to the first
$n$ iterates of $y$. It is not difficult to see that we can consider
this periodic point to be of index $1$ and such that its stable
manifold intersects the unstable manifold of $y$. So, $p \in \overline{ W^u(H)}$
and using Lyapunov stability of $H$  we know
$\overline{W^u(H)} \en H$ (see \cite{CMP} Lemma 2.1), so $p \in H$.

Since $\gamma_1$ was arbitrary, we can choose it to satisfy
$\gamma_1>\lambda$ where $\lambda$ is as in Lemma
\ref{normadiferencial}. Also, we can choose $\eps$ small so that
$\|Df^{-n}_{/F(p)}\|> \lambda^n$ contradicting Lemma
\ref{normadiferencial}.

Now, we shall study what happens if Liao's shifting Lemma can be
applied. That is, there exists $\gamma<\gamma_1<\gamma_2<1$ such
that for all $x \in H$ satisfying 

\[ \|Df^{-n}_{/F(x)}\| \geq \gamma_2^n \qquad \forall n \geq 0 \]

\noindent there exists $ y \in \omega(x)$ such that

\[ \|Df^{-n}_{/F(x)}\| \leq \gamma_1^n \qquad \forall n \geq 0 \]

So, using the Shifting Lemma we have that for every
$\gamma_2<\gamma_3< \gamma_4<1$ we have a periodic orbit $p_U$ of
$f$ contained in any neighborhood $U$ of $\Lambda$ and satisfying
that

\[ \|Df^{-n}_{/F(p)}\|\leq \gamma_4^n \]
\[ \|Df^{-n}_{/F(f^i(p))}\| \geq \gamma_3^n \]

\noindent for some $i \in 0, \ldots, \pi(p)$ (remember that $F$ is
one dimensional, so the product of norms is the norm of the
product). But since this periodic points are not very contracting in
the direction $F$, if we choose $\gamma_3 > \lambda$ (as before)
and $U$ sufficiently small to ensure that the stable manifold of
some periodic point will intersect the unstable one of a point in $H$
we reach the same contradiction as before.

\lqqd

\end{document}